\documentclass[11pt]{amsart}  
\usepackage{amssymb}
\DeclareMathOperator{\ch}{char}

\newtheorem{theorem}{Theorem}[section]

\newtheorem{lemma}[theorem]{Lemma}

\newtheorem{proposition}[theorem]{Proposition}

\newtheorem*{uncorollary}{Corollary}

\theoremstyle{definition}

\newtheorem{example}[theorem]{Example}

\theoremstyle{remark}

\numberwithin{equation}{section}

\makeatletter
\def\romenumi{%
  \def\theenumi{\roman{enumi}}%
  \def\p@enumi{\theenumi}%
  \def\labelenumi{(\@roman\c@enumi)}}
\makeatother

\newif\ifShowLabels
\ShowLabelstrue
\newdimen\mgheight
\def\marginnotes#1{%
                 \leavevmode\vadjust{\setbox0=\hbox{{\tt
                               \quad\quad {\small\textup #1}}}%
                 \mgheight=\ht0
                 \advance\mgheight by \dp0
                 \advance\mgheight by \lineskip
                 \kern -\mgheight \vbox to
                 \mgheight{\rightline{\rlap{\box0}} \vss}}}

\newcommand{\noproof}{\hfill\qedsymbol}

\begin{document}

\title[Algebras with nonmatrix and semigroup identities]{Group algebras and enveloping
algebras with nonmatrix and semigroup identities}

% author one information
\author{David~M.~Riley}  
\address{Department of Mathematics,  The University of Alabama, Tuscaloosa AL  
35487-0350, USA} 
\curraddr{} 
\email{driley@gp.as.ua.edu}
\urladdr{}
\thanks{The first author received support from NSF-EPSCoR in
Alabama and the University of Alabama Research Advisory Committee.}

% author two information
\author{Mark~C.~Wilson}  
\address{Department of Mathematics, University of Auckland, Private Bag 92019  
Auckland, New Zealand} 
\curraddr{} 
\email{wilson@math.auckland.ac.nz}
\urladdr{http://www.math.auckland.ac.nz/\~{}wilson}  
\thanks{The second author was supported by a NZST Postdoctoral Fellowship}

% 17B60, 20M07.}

\date{\today}

% at present the "communicated by" line appears only in ERA and PROC
%\commby{Ken Goodearl}

%\dedicatory{}

\begin{abstract}  
Let $K$ be a field of characteristic $p>0$.  Denote by $\omega(R)$ the
augmentation ideal of either a group algebra $R=K[G]$ or a restricted enveloping 
algebra $R=u(L)$ over $K$.   We first characterize those $R$ for which
$\omega(R)$ satisfies a polynomial identity
not satisfied by the algebra of all $2\times 2$ matrices over $K$.
Then, we examine
those $R$ for which $\omega(R)$ satisfies a semigroup identity (that is, 
a polynomial identity which can be written as the difference of two 
monomials).  
\end{abstract}  
\maketitle 

\section{Introduction and statement of results}
The structure of group algebras and restricted enveloping algebras, 
over a field $K$ of characteristic $p>0$, that are Lie nilpotent, 
Lie solvable ($p\neq2$) or satisfy the Engel condition has been 
completely determined in \cite{passi-passman-sehgal, riley-shalev:lie_structure,
sehgal:topics}.   Integral to the proof of these results was the fact that each of these 
conditions corresponds to a particular {\em nonmatrix identity},  which is to say, 
a polynomial identity not also satisfied by the algebra $M_2(K)$ of all $2\times 2$
matrices.

Nonmatrix identities in characteristic zero were studied  by Kemer
in \cite{kemer:nonmatrix},
 who showed, in particular, that nonmatrix varieties are Lie solvable. The situation in
positive characteristic is more complicated. More recently, the existence of an arbitrary
nonmatrix identity was shown in \cite{billig-riley-tasic} to be intimately related to the
existence of a group identity in the group of units of a certain class of associative 
algebras, which includes group algebras of periodic groups and  restricted enveloping
algebras of $p$-nil restricted Lie algebras. 

In this first half of this paper we completely describe enveloping algebras and group
algebras which satisfy a nonmatrix identity. In the results below, $K[G]$ and
$u(L)$ denote, respectively, the group algebra of a group $G$  and the
restricted universal enveloping algebra of the restricted Lie algebra $L$. In
addition, $\omega(R)$ represents the augmentation ideal of
$R=u(L)$ or $K[G]$,  while $\gamma^2(R)$ denotes its (associative) commutator ideal.

\begin{theorem}
\label{nonmatrix-u(L)}
The following statements are equivalent for a restricted Lie
algebra $L$ over a field of characteristic $p>0\mathrm{:}$

\begin{enumerate}
\romenumi

\item $u(L)$ satisfies a nonmatrix identity\textup{;} 
\item $\omega(u(L))$ satisfies a nonmatrix identity\textup{;}
\item $L$ contains a restricted ideal $A$ such that both $L/A$ and $A'$
 are finite-dimensional, and  $L'$ is $p$-nil of bounded index\textup{;}
\item $\gamma^2(u(L))$ is nil of bounded index.

\end{enumerate}
\end{theorem}

\begin{theorem}
\label{nonmatrix-K[G]}
The following statements are equivalent for a group $G$ and
a field $K$ of characteristic $p>0\mathrm{:}$

\begin{enumerate}
\romenumi

\item $K[G]$ satisfies a nonmatrix identity\textup{;} 
\item $\omega(K[G])$ satisfies a nonmatrix identity\textup{;}
\item $G$ has a normal subgroup $A$ such that both $A'$ and $G/A$ are
finite, and $G'$ is a $p$-group of finite exponent\textup{;}
\item $\gamma^2(K[G])$ is nil of bounded index.
\end{enumerate}

\end{theorem}  

It was conjectured in \cite{riley:boundednil} and \cite{riley:tensor},
respectively, that if $R$ and $S$ are arbitrary algebras over a 
field of characteristic $p>0$ satisfying nonmatrix identities, then 
$\gamma^2(R)$ is nil of bounded index  and the tensor product $R\otimes_K S$ 
also satisfies some nonmatrix identity. The corresponding results were shown 
to fail in characteristic 0.  Theorems~\ref{nonmatrix-u(L)} and~\ref{nonmatrix-K[G]}
verify both these conjectures for  the class of restricted enveloping algebras,
since $u(L)\otimes_Ku(M)\cong u(L\times M)$, and the class of group algebras,  since
$K[G]\otimes_KK[H]\cong K[G\times H]$. 

It will become apparent in the next section (see the paragraph after
Proposition~\ref{nonmatrix}) that if $R$ is either an  ordinary enveloping algebra over an
arbitrary field, or a group algebra over a field of  characteristic zero, and $R$
satisfies a nonmatrix identity, then
$R$ must be commutative. 

In the second half of this article we look at  restricted
enveloping algebras and group algebras which satisfy a semigroup identity.  
A {\em semigroup identity} is a polynomial identity of the form 
$$f(x_1,\dots,x_n)=w_1(x_1,\dots,x_n)-w_2(x_1,\dots,x_n)=0,$$
where $w_1$ and $w_2$ are monomials.
It follows from well-known results that  a semigroup
identity is a nonmatrix identity precisely  when the base field $K$ is infinite. In fact,
by a result of Jones \cite{jones:varieties}, every proper variety of groups can contain
only finitely many simple groups. Hence the general linear group $GL_2(K)$ does not
satisfy any group identity unless $K$ is finite. This last fact is also a special case of
\cite[Theorem 1.2]{billig-riley-tasic}. The present authors studied arbitrary algebras
$R$ over an infinite field that satisfy a semigroup identity in
\cite{riley-wilson:semigroup}.   In particular, we proved that if $R$ is unital then $R$
satisfies a semigroup  identity if and only if $R$ satisfies the Engel condition:
$$[x, _my]=[x,\underbrace{y,\dots,y}_m]=0,$$
for some $m$.
Hence, using the characterization of the Engel condition in restricted enveloping 
algebras and group algebras in \cite{riley-shalev:lie_structure} and \cite{sehgal:topics},
respectively, we immediately obtain the equivalence of the first three conditions 
in each of the following two theorems:

\begin{theorem}
\label{semigroup-u(L)}  
The following statements are equivalent for a restricted Lie algebra $L$ 
over an infinite field of characteristic $p>0\mathrm{:}$
\begin{enumerate}
\romenumi

\item $u(L)$ satisfies a semigroup identity\textup{;}
\item $u(L)$ satisfies the Engel condition\textup{;} 
\item $L$ is nilpotent, $L'$ is $p$-nil of bounded
index, and $L$ contains a restricted ideal $A$ 
such that both $L/A$ and $A'$ are finite-dimensional\textup{;}
\item $L'$ and $L/Z(L)$ are $p$-nil of bounded index
and $L$ contains a restricted ideal $A$ such that both $L/A$ and
$L'$  are finite-dimensional.

\end{enumerate} 
\end{theorem}

\begin{theorem}
\label{semigroup-K[G]}  
The following statements are equivalent for a group $G$ 
and an infinite field $K$ of characteristic $p>0\mathrm{:}$
\begin{enumerate}
\romenumi

\item $K[G]$ satisfies a semigroup identity\textup{;}
\item $K[G]$ satisfies the Engel condition\textup{;}
\item $G$ is nilpotent and contains a normal subgroup $A$ such that both
$G/A$ and $A'$ are finite $p$-groups\textup{;} 
\item $G'$ and $G/Z(G)$ are $p$-groups of finite exponent and
$G$ contains a normal subgroup $A$ such that both $G/A$ and $A'$ are finite.

\end{enumerate} 
\end{theorem}

Note that since every finite ring satisfies an identity of the form 
$x^k-x^l=0$, these results do not hold when $K$ is finite.
Observe also that the corresponding conditions (iv) provide an exact 
analogue between the group algebra and restricted enveloping algebra cases.
These alternate characterizations will prove useful below.
The proof of their equivalence requires some work and will be given later.

For nonunital algebras, semigroup identities are more difficult to 
classify. In particular, the existence of a semigroup identity does not generally imply
the Engel condition. In fact, this implication does not even hold for augmentation
ideals: 

\begin{example}
\label{e:2-D}
 Let $K$ be a field of characteristic 2, and let $L$ 
be the 2-dimensional restricted Lie algebra over $K$ generated by
$a,b$ with relations $[a,b]=a$, $a^2=0$, $b^2=b$.
Then $\omega(u(L))$ satisfies the semigroup identity 
$$wxyz-wyxz=0,$$
despite the fact that $\omega(u(L))$ does not satisfy the Engel
condition.
\end{example}

\begin{example}
\label{e:dihedral} Let $K$ be a field of characteristic 3, and let $G$ be
the dihedral group of order 6.  Then $\omega(K[G])$ satisfies the
semigroup identity 
$$yxy^2-y^2xy=0,$$
and yet $\omega(K[G])$ does not satisfy the Engel condition.
\end{example}

We shall see in Section~\ref{s:semigroup} though, that these two examples are in some
sense canonical.

\begin{theorem}
\label{semigroup-omega(u(L))}

Let $L$ be a restricted Lie algebra over an infinite field 
of characteristic $p>2$.
Then the following statements are equivalent\textup{:}
\begin{enumerate}
\romenumi

\item $\omega(u(L))$ satisfies a semigroup identity\textup{;}
\item $\omega(u(L))$ satisfies an identity of the form $y^m[x, _my]y^m=0$\textup{;}
\item $L'$ and $L/Z(L)$ are both $p$-nil of bounded index and 
$L$ contains a restricted ideal $A$ of finite codimension such that $L'$ is 
finite-dimensional\textup{;} 
\item $u(L)$ satisfies the Engel condition.

\end{enumerate}
\end{theorem}

\begin{theorem}
\label{semigroup-omega(K[G])}
Let $K$ be an infinite field of characteristic $p>0$, and suppose that $p=2$ or $G$ is a
group with no $2$-torsion. Then the following statements are equivalent\textup{:}
\begin{enumerate}
\romenumi

\item $\omega(K[G])$ satisfies a semigroup identity\textup{;}
\item $\omega(K[G])$ satisfies an identity of the form $y^m[x, _my]y^m=0$\textup{;}
\item $G'$ and $G/Z(G)$ are $p$-groups of finite exponent and $G$ contains a 
normal subgroup $A$ of finite index such that $A'$ is finite\textup{;}
\item $K[G]$ satisfies the Engel condition.

\end{enumerate}
\end{theorem}

\section{Nonmatrix identities}
In this section we prove Theorems~\ref{nonmatrix-u(L)} and~\ref{nonmatrix-K[G]}.
Throughout, $K$ is a field, $G$ a group and $L$ a restricted Lie
$K$-algebra. For each (ordinary) subalgebra $H$ of $L$, $H_p$ denotes its $p$-hull
inside $L$.

We begin by collecting some important facts about algebras  satisfying a
nonmatrix identity. 

\begin{proposition}
\label{nonmatrix}
Let $R$ be an algebra over a field $K$. 
If $R$ satisfies a nonmatrix identity then the following statements hold\textup{:}
\begin{enumerate}
\romenumi
\item The unital hull of $R$ also satisfies a nonmatrix identity.

\item $\gamma^2(R)$ is nil, and is nilpotent if $R$ is finitely
generated over $K$.

\item If $\ch K=p>0$ then $R$ satisfies an identity of the form
$([x,y]z)^{p^t}=0$, so that $\gamma^2(R)$ is generated by nilpotent 
elements of bounded index.

\item If $\ch K=p>0$ and  $R=\omega(K[G])$ then $G'$
is a $p$-group generated by elements of bounded $p$-power order.

\item If $\ch K=p>0$ and $R=\omega(u(L))$ then $(L')_p$ is a $p$-nil
restricted Lie algebra generated by $p$-nilpotent elements of bounded index.  
\end{enumerate}

\end{proposition}

\begin{proof}

Let $S$ be the unital hull of $R$, that is, $S=R$ if $R$ is unital and $S=R\oplus
K\cdot 1$ otherwise. Then $R$ is an associative ideal of $S$ such that$[S,S]S\subseteq R$.
 Let $f(x_1,\dots ,x_n)=0$ be a nonmatrix identity satisfied by $R$. It follows that $S$
satisfies the identity 
$$f([x_1,y_1]z_1,\dots ,[x_n,y_n]z_n)=0.$$
Now this last identity is also nonmatrix, since 
$$a=[e_{12},e_{21}](e_{11}-e_{22})a$$
for every $a\in M_2(K)$, and this proves (i).

Henceforth, we assume that $R$ is unital.

Let $f$ be the given nonmatrix identity, and $J(R)$ the
Jacobson radical of $R$. We first show that $R/J(R)$ is commutative. Now
$R/J(R)$ is a semiprimitive PI-algebra, and hence a subdirect product of
primitive PI-algebras, each of which satisfies $f$.  By a theorem of
Kaplansky, a primitive PI-algebra is simple and finite-dimensional over
its centre, and thus has the form $M_n(D)$ where $D$ is a division algebra
finite-dimensional over its centre $F$. Since $M_2(F)$ does not satisfy
$f$, $n=1$. If $D$ is finite then $D=F$, while if $D$ is infinite then
$D$ satisfies the same identities as some $M_m(F)$. The nonmatrix
condition gives $m=1$ and so either way $D=F$. Thus $R/J(R)$ is a
subdirect product of fields and hence commutative. It
follows that $\gamma^2(R)\subseteq J(R)$.

By a theorem of Braun \cite{braun:nilpotency}, if $A$ is a finitely
generated PI-algebra then $J(A)$ is nilpotent. Let $A$ be the relatively
free algebra on 3 generators in the variety defined by the given
nonmatrix identity $f$.  Then  $\gamma^2(A)\subseteq J(A)$ is nilpotent
and so $A$ satisfies an identity of the form $([x,y]z)^{m}=0$. It
follows that $R$ also satisfies this same identity. This yields (ii) and
(iii). If $R=u(L)$ then $L'\subseteq \gamma^2(R)$, whereas if $R=K[G]$ then
$G'\subseteq 1+\gamma^2(R)$, yielding (iv) and (v).

\end{proof}

If an ordinary enveloping algebra $U(L)$ satisfies a nonmatrix identity
then $L'$ is nil. By the Poincar\'{e}-Birkhoff-Witt (PBW) theorem, $L'=0$
and $L$ is abelian. For similar reasons, $K[G]$ satisfies a nonmatrix
identity in characteristic $0$ if and only if $G$ is abelian. 
Thus we shall focus on modular group algebras and restricted enveloping algebras from now
on. Because of the strong parallels between the group algebra and restricted
enveloping algebra cases, we prove only the restricted enveloping algebra 
case in detail.  

For the rest of this section, $K$ will be assumed to be of characteristic $p>0$.
\subsection*{Enveloping algebras}

We adapt a lemma of Passman \cite[Lemma 3.2]{passman:units2} to
the enveloping algebra situation.

\begin{lemma}
\label{expand}
Let $A$ be a restricted ideal of an abelian restricted Lie
algebra $L$ such that $\dim L/A$ is finite. If $I$ is an $L$-stable
ideal of $u(A)$ then the ideal $J=Iu(L)=u(L)I$ of $u(L)$
generated by $I$ is  nil of bounded index if (and only if) $I$ is.

\end{lemma}

\begin{proof} Let $n=\dim L/A$. Jacobson's theorem (the analogue of the
PBW theorem for restricted Lie algebras) implies
that as a left $u(A)$-module, $u(L)$ is free of rank $q=p^n$.
Indeed, if $\{x_1,\dots ,x_n\}$ is an ordered basis for a vector space
complement for $A$ in $L$, the monomials $x_1^{i_1}\cdots x_n^{i_n}$
form a basis for $u(L)$ over $u(A)$. 
 
The right regular representation of $u(L)$ embeds $u(L)$ into $M_q(u(A))$
in such a way that for all  basis monomials $\nu$, $\nu
x=\sum_{\mu}x_{\nu\mu}\mu$ for $x\in u(L)$. If
$x\in J$ then also $\nu x\in J$.  Since $J=\bigoplus_{\mu} I\mu$ we have
$\sum_{\mu}x_{\nu\mu}\mu\in \bigoplus_{\mu} I\mu$.  Since all
$x_{\nu\mu}\in u(A)$, freeness yields
$x_{\nu\mu}\in I$ and so $J$ embeds in $M_q(I)$.

If $I$ is nil of bounded degree $p^t$, then since every element of
$M_q(I)$ is algebraic of degree at most $q$ over the central 
subalgebra $I$, it follows easily that $M_q(I)$ is nil of bounded
index at most $qp^t$.

\end{proof}

\begin{proof}[Proof of Theorem~\ref{nonmatrix-u(L)}]
Clearly (i) is equivalent to (ii) by Proposition~\ref{nonmatrix}, while the fact that (iv)
implies (i) is immediate.  

Recall that for an (ordinary) ideal $J$ of $L$, the ideal generated by 
$J$ in $u(L)$ is $\omega(u(J_p))u(L)=u(L)\omega(u(J_p))$.
By Jacobson's theorem, this ideal is nilpotent if and only if 
$J_p$ is finite-dimensional and $p$-nil.

Results of Passman \cite{passman:u(L)-PI} and Petrogradski \cite{petrogradski}, 
yield that $u(L)$ satisfies a polynomial identity if and only if 
$L$ contains a restricted ideal $A$ such that $\dim L/A<\infty$,
$\dim A'<\infty$ and $(A')_p$ is $p$-nil. 
Suppose that condition (i) holds, so that $L$ has a restricted ideal 
$A$ with these basic properties.  In particular, $A$ is solvable
since by Engel's theorem it is nilpotent-by-abelian.
To show (i) implies (iii), it remains to show that $(L')_p$ is $p$-nil of bounded
index. Certainly $(L')_p$ is $p$-nil and generated by $p$-nilpotent elements of 
bounded index by Proposition~\ref{nonmatrix}. 
Thus $(L')_p/(A\cap (L')_p)$ is  finite-dimensional and $p$-nil, 
so that $(L')_p$ is solvable. Consequently, $L$ is solvable.   
Now arguing by induction on the derived length of $L$ enables
us to reduce to the case when $L'$ is abelian. 
But in this case $(L')_p$ is clearly $p$-nil of bounded index,
since it is generated by elements with that property.

Finally, suppose that (iii) holds. We show that $\gamma^2(u(L))$ is nil
of bounded index. We make a series of reductions based on the
fact that the class of algebras which are nil of bounded index is closed
under extensions.  
Because $A'$ generates a nilpotent ideal in $u(L)$, it suffices for us to 
assume $A$ is abelian. Now put $B=[A,L]$. Then $B$ is an ideal of $L$
contained in $A\cap L'$ and so $B_p$ is $p$-nil of bounded index. 
Since $A$ is abelian, it follows that $B$ generates an ideal $I$ of 
$u(A)$ which is $L$-stable and nil of bounded index. By Lemma~\ref{expand} the same
is true of the expanded ideal $Iu(L)$ of $u(L)$. Thus we may assume that $B=0$ and
hence that $A$ is central in $L$. This implies that $L'$ is finite-dimensional
and so $L'$ generates a nilpotent ideal of $u(L)$. Finally, we can assume
that $L$ is abelian and the result follows immediately.
\end{proof}    

\begin{uncorollary}
$u(L)$ satisfies a nonmatrix identity if and only if it satisfies some polynomial identity
and $L'$ is $p$-nil of bounded index.
\noproof
\end{uncorollary}

\subsection*{Group algebras}

\begin{proof}[Proof of Theorem~\ref{nonmatrix-K[G]}]
The proof is entirely analogous to the restricted enveloping algebra case, and is in
fact simpler since there  is no need to trouble ourselves with any analogue of the
$p$-hull. 

The main ingredients are: 
\begin{itemize}
\item By a result of Passman \cite[Corollary 5.3.10]{passman:asgr}, $K[G]$
satisfies a polynomial identity if and only if $G$ has a subgroup
$A$ such that $G/A$ is finite and $A'$ is a finite $p$-group.
\item A normal subgroup $H$ of $G$ generates a nilpotent ideal of
$K[G]$ if and only if $H$ is a finite $p$-group. 
\item In place of Lemma~\ref{expand}, we use the group algebra
analogue \cite[Lemma 3.2]{passman:units2}.
\end{itemize}
\end{proof}

\begin{uncorollary}
$K[G]$ satisfies a nonmatrix identity if and only if it satisfies some polynomial identity
and $G'$ is a $p$-group of finite exponent.
\noproof
\end{uncorollary}

\section{Semigroup identities}
\label{s:semigroup}
In this section we prove Theorems~\ref{semigroup-u(L)}, \ref{semigroup-K[G]}, 
\ref{semigroup-omega(u(L))} and \ref{semigroup-omega(K[G])}, and justify
Examples~\ref{e:2-D} and \ref{e:dihedral}. The same notational hypotheses as in the
previous section are in force throughout. In addition, the field $K$ is assumed to be
infinite of characteristic $p>0$.

\subsection*{Enveloping algebras}

\begin{proof}[Proof of Theorem~\ref{semigroup-u(L)}]
It remains to prove the equivalence of condition (iv) and condition (iii), say.
That (iii) implies (iv) is obvious.  Conversely, assume that
(iv) holds.  Replacing $A$ by the centralizer of $A'$ in $A$,
we may assume that $A$ is a nilpotent restricted
ideal of $L$ containing $Z(L)$.  But then $L$ is 
nilpotent-by-(finite-dimensional and $p$-nil), and hence
nilpotent by \cite{shalev:PI-grK[G]}.
\end{proof} 

\begin{proposition}
\label{2-D} Let $p=2$ and let $D$ be the $2$-dimensional restricted 
Lie algebra over $K$ generated by
$a,b$ with relations $[a,b]=a$, $a^2=0$, $b^2=b$.
Then $\omega(u(D))$ satisfies the semigroup identity 
$$wxyz-wyxz=0.$$
\end{proposition}

\begin{proof} By Jacobson's theorem, $\omega(u(D))$ has $K$-basis
$$\{a,b,ab\}.$$
Because the given semigroup identity is multilinear, it suffices to
show that it holds on basis elements.
Observe
$$[a,b]=a,\quad [a,ab]=a^2=0,\quad [ab,b]=ab=(1+b)a,$$
each of which is annihilated by $a$ on the right and left.
The remaining possibilities are:
$$bab=b(1+b)a=2ba=0,\quad b[(1+b)a]b=2bab=0.$$
\end{proof}

\begin{proposition} If $\omega(u(L))$ satisfies an identity of the form 
$y^m[x, _m y]y^m=0$ then either $u(L)$ satisfies the Engel condition or 
$p=2$ and $L$ contains a restricted subalgebra isomorphic to $D$.
\end{proposition}
\begin{proof} Since $y^m[x, _m y]y^m=0$ is a nonmatrix identity,
it is clear from Theorems~\ref{nonmatrix-u(L)} and \ref{semigroup-u(L)} that
$u(L)$ satisfies the Engel condition as soon as $L/Z(L)$ is
$p$-nil of bounded exponent.   Because
$$y^m[x, _{m+1}y]y^m=[y^m[x, _m y]y^m,y]=0$$
is also an identity for $\omega(u(L))$, we may
assume that $m=p^t$ for some $t$.
Then $\omega(u(L))$ satisfies the identity   
$$y^{p^t}[x,y^{p^t}]y^{p^t}=y^{p^t}[x, _{p^t} y]y^{p^t}=0.$$
Suppose now that $a,c$ are arbitrary elements in $L$ and put
$b=c^{p^t}$.
Then
$$[a,b]b^2-[a,b,b]b=b[a,b]b=0.$$
Assume, for the moment, that $p>2$.
Then it follows easily using Jacobson's theorem that $[a,b,b]=0$\textup{;}
hence, $L^{p^{t+1}}\subseteq Z(L)$, as required.

Assume, finally, that $p=2$ and, without loss of generality, 
$[a,_3 b]\neq0$ for some $a,c$ as above. 
Then, by Jacobson's theorem, there exist 
$\lambda,\mu,\nu\in K$ such that
$[a,b,b]=\lambda [a,b]+\mu b+\nu b^2.$
It follows that $[a,b,b,b]=\lambda [a,b,b]$, where $\lambda\neq0$.
Replacing $a$ by $[a,b,b]$ and subsequently $b$ by $\lambda^{-1} b$,
we can assume $[a,b]=a$.  Therefore 
$$a(b^2-b)=[a,b]b^2-[a,b,b]b=0,$$
so that $b^2=b$.
Clearly now, $a$ and $b$ generate a restricted subalgebra of 
$L$ isomorphic to $D$, as required.
\end{proof}

\begin{proof}[Proof of Theorem~\ref{semigroup-omega(u(L))}]
According to
\cite{riley-wilson:semigroup}, if an algebra over an infinite field
satisfies a semigroup identity then it satisfies an identity
of the form $$y^m[x,_m y]y^m=0,$$
for some $m$.  Now combining the preceding
proposition with Theorem~\ref{semigroup-u(L)} readily yields the result.
\end{proof}

\subsection*{Group algebras}

\begin{proof}[Proof of Theorem~\ref{semigroup-K[G]}]
It remains to prove the
equivalence of condition (iv) and conditions (i)-(iii).
That  condition (ii) implies (iv) is clear from Theorem~\ref{nonmatrix-K[G]}
since if $g,h\in G$ then 
$$gh^{p^t}-h^{p^t}g=[g,h^{p^t}]=[g, _{p^t}h].$$  
Assume now that only (iv) holds.  Replacing $A$ by the centralizer 
of $A'$ in $A$, we may assume that $A$ is a nilpotent normal
subgroup of $G$ containing $Z(G)$.  But then $G/Z(G)$ 
is  a $p$-group of finite exponent that is an extension of a 
nilpotent group by a finite $p$-group.  
Such a group is nilpotent by \cite{baumslag}.  Hence,
$G$ itself is nilpotent and (iii) holds.
\end{proof}

\begin{proposition}
\label{dihedral}
Let $p=3$ and let $G$ be the dihedral group of order $6$. Then $\omega(K[G])$ satisfies
the identity $yxy^2=y^2xy$.
\end{proposition}

\begin{proof}
Write $R=\omega(K[G])$. We use the presentation $G=\langle a,b | a^2=b^3=aba=1\rangle$.
Let $d=b-b^{-1}\in R$ and let $e=(1-a)/2\in R$. It is easy to see that $e^2=e, d^3=0$ and 
$d^2$ is central in $R$. Furthermore, $ed=d(1-e)$. A $K$-basis for $R$ is $\{e,d^2,
ed^2,d,ed\}$. Let $R_0$ be the span of $e, d ^2, ed^2$ and let $R_1$ be the span of $d,
ed$. Then $R=R_0\dotplus R_1$ is a superalgebra and the following basic properties hold:
\begin{itemize}
\item $R_0$ is commutative;
\item $(R_1)^3=0$;
\item $R_0R_1R_0=0$.
\end{itemize}

Consider the given identity in the form $y[x,y]y=0$. Since the identity is linear in $x$
it suffices to assume that $x\in R_0$ or $x\in R_1$. Fix $y=s_0+s_1\in R_0\dotplus R_1$.

Suppose that $0\neq x=r_0\in R_0$. Then using the above properties we obtain
$$
y[x,y]y=s_0[r_0,s_1]s_1+s_1[r_0,s_1]s_0.
$$
Since $d^2$ kills all basis elements except $e$, and the last expression is homogeneous
in $r_0$ and $s_0$, without loss of generality we may assume that $r_0=s_0=e$. Then the
above simplifies to 
$$
y[x,y]y=e[e,s_1]s_1+s_1[e,s_1]e=es_1^2-s_1^2e=0.
$$

Suppose that $0\neq x=r_1\in R_1$. Then the basic properties above yield
$$
y[x,y]y=[r_1,s_1]s_0^2+s_0[r_1,s_0]s_1+s_1[r_1,s_0]s_0.
$$ 

Since $x$ has degree 1 in $d$ we may again assume that $s_0=e$, which yields
$$
y[x,y]y=[r_1,s_1]e+e[r_1,e]s_1+s_1[r_1,e]e=r_1s_1e-s_1r_1e-er_1s_1+s_1r_1e=0.
$$
\end{proof}

\begin{proposition} If $\omega(K[G])$ satisfies an identity of the form 
$y^m[x,_m y]y^m=0$ then either $K[G]$ satisfies the Engel condition or $p>2$ and $G$
contains a subgroup isomorphic to the dihedral group $D_{2p}$ of order $2p$.
\end{proposition}
\begin{proof} As above, if $G/Z(G)$ is a $p$-group of bounded
exponent, then the assumption on $\omega(K[G])$ forces $K[G]$ to satisfy
the Engel condition.

Assume $m=p^t$.  We claim that given $h\in G$, either $h^{p^{t+1}}\in Z(G)$ or
$h^{2p^t}=1$.  Let $g,h\in G$.  Substituting $x=1-g,y=1-h$ into the  polynomial
identity yields
$$(1-h^{p^t})[h^{p^t},g](1-h^{p^t})=(1-h)^{p^t}[g,(1-h)^{p^t}](1-h)^{p^t}=0.$$
Expanding, we obtain
$$h^{p^t}g-gh^{p^t}-h^{2p^t}g+gh^{2p^t}+h^{2p^t}gh^{p^t}-h^{p^t}gh^{2p^t}=0.$$

Consider first the case when $p>5$.  Then by the linear independence of group
elements we have:
$$\{h^{p^t}g,gh^{2p^t},h^{2p^t}gh^{p^t}\}=\{gh^{p^t},h^{2p^t}g,h^{p^t}gh^{2p^t}\}.$$
By inspection, it follows that either $h^{p^t}$ commutes with $g$ or $h^{2p^t}=1$,
as claimed.

Next assume $p=3$.  One additional possibility can occur:
$$h^{3^t}g=gh^{2\cdot3^t}=h^{2\cdot3^t}gh^{3^t},
\quad gh^{3^t}=h^{2\cdot3^t}g=h^{3^t}gh^{2\cdot3^t}.$$
But then $h^{3^t}g=h^{2\cdot3^t}gh^{3^t}$ implies $g=h^{3^t}gh^{3^t}$, so that
$g=(h^{3^t}g)h^{3^t}=(gh^{2\cdot3^t})h^{3^t}$ implies that $h^{3^{t+1}}=1$.

Finally, in the $p=2$ case the fact that the group element terms above must cancel 
pairwise eventually shows that $h^{2^{t+1}}$ commutes with $g$, thus establishing 
the claim for all characteristics.  This also completes the proof in the
case $p=2$\textup{;} henceforth, we assume $p>2$.

Now let $P$ denote the set of $p$-elements in $G$.  Since $\omega(K[G])$ 
satisfies a nonmatrix identity, $G'$ is a $p$-group by Theorem~\ref{nonmatrix-K[G]}.
It follows immediately that $P$ is a normal subgroup of $G$.

According to the claim above, either $G/Z(G)$ is a $p$-group of 
finite exponent or $G$ contains an element $\sigma$ of order $2$.  
Consider the case that $\sigma$ centralizes $P$.  
We claim that $\sigma$ must be central.  
Indeed, let $\tau\in G$. Then since $G'$ is a $p$-group, 
$\sigma^\tau=\sigma g$ for some $p$-element $g$.  Thus
$$g^2=\sigma^2 g^2=(\sigma g)^2=(\sigma^\tau)^2=1,$$
so that $g=1$ since $p$ is odd.

Now consider the case when $\sigma$ does not centralize $P$.
Then there is some $p$-element $g$ in $G$ such that
$g^\sigma\neq g$.  Put $h=g^{-1}g^\sigma$, another
$p$-element in $G$.  Then
$$h^\sigma h=(g^{-1}g^\sigma)^\sigma g^{-1}g^\sigma=(g^\sigma)^{-1}gg^{-1}g^\sigma=1.$$ 
If the order of $h$ is $p^n$ and $n\ge2$, then replacing
$h$ by $h^{p^{n-1}}$ we may assume $h$ has order exactly 
$p$ and $h^\sigma=h^{-1}$.  Consequently, the subgroup 
$G$ generated by $\sigma$ and $h$ is isomorphic to $D_{2p}$.
\end{proof}

\begin{proof}[Proof of Theorem~\ref{semigroup-omega(K[G])}]
According to
\cite{riley-wilson:semigroup}, if an algebra over an infinite field
satisfies a semigroup identity then it satisfies an identity
of the form $$y^m[x,_m y]y^m=0,$$ for some $m$.  Now combining the preceding proposition
with Theorem~\ref{semigroup-K[G]} readily yields the result.
\end{proof}  

\section{Further comments}
The complete analogy between Theorems~\ref{nonmatrix-u(L)} and \ref{nonmatrix-K[G]}
suggests that a more general result may be possible, perhaps in some class of Hopf
algebras containing both examples.

We strongly suspect that $\omega(K[D_{2p}])$ satisfies a semigroup identity for fields
$K$ of arbitrary prime characteristic $p$, but do not have a general proof. The
difficulties encountered in even this seemingly simple case lead us to believe that the
full characterization of augmentation ideals satisfying a semigroup identity will require
considerable additional effort.

\bibliographystyle{amsalpha}  
\bibliography{biblio}  
  
\end{document}